\documentclass{radhazu}
\usepackage{radhazu}
\usepackage{graphicx}

\newcommand{\Z}{\mathbb{Z}}
\newcommand{\Q}{\mathbb{Q}}

\begin{document}

\title[An elliptic curve over $\mathbb{Q}(u)$ with torsion $\mathbb{Z}/ 4\mathbb{Z}$ and rank $6$]{An elliptic curve over $\mathbb{Q}(u)$ with torsion $\mathbb{Z}/ 4\mathbb{Z}$ and rank $6$}
\author[A. Dujella]{Andrej Dujella}
\address{Department of Mathematics \\ Faculty of Science \\ University of Zagreb \\ Bijeni\v{c}ka cesta 30, 10000 Zagreb, Croatia}
\email[A. Dujella]{duje@math.hr}
\author[J. C. Peral]{Juan Carlos Peral}
\address{Departamento de Matem\'aticas \\ Universidad del Pa\'is Vasco \\ Aptdo. 644, 48080 Bilbao, Spain}
\email[J. C. Peral]{juancarlos.peral@ehu.es}
\thanks{     }

\abstract{
In this paper, we present the construction of
an elliptic curve over $\mathbb{Q}(u)$ with torsion group $\mathbb{Z}/ 4\mathbb{Z}$ and rank $6$.
Previously only rank $5$ examples for such curves were known.}

\maketitle

\begin{center}
{\it Dedicated to Marko Tadi\'c on the occasion of his $70$th birthday}
\end{center}

\section{Introduction}

Let $E$ be an elliptic curve over $ \mathbb{Q}$.
By {the Mordell-Weil theorem}, the group $E( \mathbb{Q})$ of rationals points
on $E$ is a finitely generated Abelian group. Hence, it has the form
$E( \mathbb{Q}) \cong E( \mathbb{Q})_{\rm tors} \times { \mathbb{Z} }^r$,
where $E( \mathbb{Q})_{\rm tors}$ is called the torsion group of $E$, and the non-negative integer $r$
is called the rank of $E$.

In 1978, Mazur \cite{Ma} proved that there are precisely $15$ possible torsion groups for elliptic curves over $\mathbb{Q}$:
\begin{eqnarray*}
\mathbb{Z}/k\mathbb{Z},& &  \mbox{for $k=1,2,3,4,5,6,7,8,9,10,12$,} \\
\mathbb{Z}/2\mathbb{Z}  \times \mathbb{Z}/k\mathbb{Z},& &  \mbox{for $k=2,4,6,8$.}
\end{eqnarray*}

On the other hand, it is not known which values of rank $r$ are
possible for elliptic curves over $ \mathbb{Q}$.
It has been conjectured that there exist elliptic curves of arbitrarily high rank,
even that for each of the torsion groups in Mazur's theorem, there exist curves with that torsion group
and arbitrarily large rank.
However, there are also recent heuristic arguments that suggest the boundedness of the
rank of elliptic curves. According to this heuristic,
only a finite number of curves would have rank higher than 21 (see \cite{PPVW}).
The current record is an
example of an elliptic curve over $\mathbb{Q}$ with rank $\geq 28$,
found by {Elkies} in 2006 (see \cite{El} and also \cite{KSW} where it is shown that,
assuming the generalized Riemann hypothesis, Elkies's curve has rank equal to $28$).

The first step in finding elliptic curves with large rank is usually
constructing a parametric family of elliptic curves,
i.e.~an elliptic curve over the field of rational function $\mathbb{Q}(u)$,
with the considered property (e.g.~with a given torsion group).
Thus, it is of interest to construct elliptic curves over  $\mathbb{Q}(u)$
with a given torsion group and relatively large generic rank.

In this paper, we will present the construction of
an elliptic curve over  $\mathbb{Q}(u)$ with torsion group $\Z / 4 \Z $ and rank $6$,
which is the current record for the rank over $\mathbb{Q}(u)$ for this torsion group.
Let us mention that the largest known rank for an elliptic curve with torsion group  $\Z / 4 \Z $
over $\mathbb{Q}$ is $13$ (see \cite{ElKl}).
Information on current rank records for all $15$ possible torsion groups can be found in \cite{Du1}.

\section{General equation of curves with torsion group $ \Z / 4\Z$}

The general model of an elliptic curve with torsion group $ \Z / 4\Z$
is given by
\[
Y^2+ a X Y+ a b Y=X^3+ b X^2,
\]
where $ab(a^2-16b)\neq 0$ (see, e.g. \cite{El}).
A torsion point of order $4$ in this model is  $(0,0)$.

We will use the construction of Elkies \cite{El} for curves with torsion group  $ \Z / 4\Z$.
Elkies noticed that this torsion and rank $4$ can be obtained for some elliptic $K3$ surfaces.
In this case, the maximum rank is obtained with the following types of reducible fibers for such a surface:
four of type $I_4$,  two of type $I_2$ and four of type $I_1$, so giving a contribution to
the N\'eron-Severi group of $4(4-1)+2(2-1))=14$,
hence the rank over this surface is at most $20- 2 - 14=4$
(since the N\'eron-Severi rank of a $K3$ surface is at most 20, see, e.g. \cite{ScSh}).
So in this sense, the Elkies example is optimal.

\section{Elkies model for rank $4$}

By the result of Inose \cite{In}, we know that the N\'eron-Severi rank can be 20 only if the corresponding
discriminant is one of the thirteen discriminants of imaginary quadratic orders of class number 1.
 Elkies \cite{El} has shown that the discriminant $-163$ surface does have an elliptic model that attains rank $4$
 with torsion group $ \Z / 4\Z$ for the following values of parameters $a$ and $b$:
   \begin{align*}
a&=(8 t-1)(32 t +7), \\
b&=8(t+1)(15 t-8)(31 t-7).
\end{align*}
  With a simple change of variables, the surface can be written as
\[
Y^2=X^3+ (a^2- 8 b) X^2+ 16 b^2 X.
\]
Inserting the values of $a$ and $b$ given above, one gets the following elliptic $K3$ surface $E$ 
   \begin{align*}
E:\,\, Y^2&=X^3+(65536 t^4- 17472 t^3- 10176 t^2+ 18672 t -3535) X^2\\
&\,\,\,\mbox{}+1024 (t+1)^2 (15 t -8)^2(31 t-7)^2 X. 
\end{align*}
It has  torsion group $ \Z / 4\Z$ and rank $4$.
A torsion point of order $4$ in this model is
 \[
(32(t+1)(15 t-8)(31 t-7), 2^5 (1+t) (-1+8 t) (-8+15 t) (-7+31 t) (7+32 t) )\]
and
the $X$-coordinates of four independent points of infinite order are 
   \begin{align*}
X_1&=-361 (t+1) (31 t-7),\\
X_2&=-4(t+1)(15 t-8)(16 t -7)^2,  \\
X_3&=-16 (t+1)(8 t+7)^2(15 t -8), \\
X_4&=4(15 t-8)(16 t+1)^2 (31 t-7).
\end{align*}

In \cite{ElKl}, it is quoted that for
 \begin{align*}
a&=4(9+80 t), \\
b&=1/2(-2+t)(-81+ 2t)(-1+ 2t)(-1+18 t).
\end{align*}
the resulting $K3$ surface also has rank $4$.

Previous to the Elkies construction in 2007, in 2004, Kihara \cite{Ki1,Ki2} and Lecacheux \cite{Le}
found curves with rank $4$ and $5$ over $\Q(u)$.
The coefficients in the Kihara and Lecacheux constructions are much larger than those in Elkies examples,
and so they are less suited for finding particular good examples of high-rank curves.
In fact, the family with rank $4$ in the Elkies paper is $y^2=x^3+ A(t) x^2+B(t) x$,
where $A(t)$  is a polynomial of degree $4$
and $B(t)$ is a polynomial of degree  $6$, and the examples with rank $5$ have polynomial coefficients
of degree $8$  and $16$, while the corresponding coefficients in the Kihara family have degrees $52$ and $102$,  respectively. In 2016, Khoshnam and Moody \cite{KhMo} gave an example of rank $5$  over $\Q(u)$ with
a simpler version of the Kihara method. In fact, the coefficients $A$ and $B$ in this curve  are polynomials
of degrees  $19$ and $38$, respectively.

The object of this note is to construct an example of an elliptic curve over $\Q(u)$
with  torsion $ \Z / 4\Z$ and rank $6$.
It will be obtained as a specialization of the Elkies surface.

\section{A  curve over $\Q(u)$ with torsion $ \Z / 4\Z$  and rank $6$}

Elkies \cite{El} mentions, without explicitly  writing such examples, that there are several quadratic sections
giving  curves over $\Q(u)$  with  rank $5$  for this torsion and several combinations of pairs of quadratic
sections leading to infinite families of curves with rank $6$ parametrized by the points of  elliptic curves of
positive rank.

In \cite{DuPe}, we made explicit the ideas of Elkies, and presented there $25$ elliptic curves over $\Q(u)$
with  rank $5$ and several infinite families of curves with rank $6$ parametrized by elliptic curves with
positive rank.

Now we have many more examples  of quadratic sections over $\Q(u)$  leading to  rank $5$ curves.
In one of these  curves, we found  a further quadratic section  leading to the curve over $\Q(u)$
having  rank $6$.

In fact, observe that imposing
$$\frac{-64 (1 + t)^2 (-4 + 7 t) (4 + 17 t)}{(1 + 4 t)^2}$$
as the $X$-coordinate of a new point in the Elkies surface $E$ is equivalent to the condition
$-(-4+7 t) (4+17 t)=\Box$, which can be solved with
$$t \mapsto \frac{4 (-1 + u^2)}{(17 + 7 u^2)}.$$
The resulting curve has rank $5$.
(For a general strategy for imposing new points which might lead to conditions equivalent to genus $0$ or $1$
curves, see, e.g. \cite{DuPe2}.)


On the other hand, imposing
$$\frac{576 (-4 + 7 t) (-8 + 15 t)^2 (-1324 + 5551 t)}{49 (-39 + 28 t)^2}$$
as the $X$-coordinate of a new point in the Elkies surface $E$ is equivalent to the condition
$(-4 + 7 t) (-1324 + 5551 t)=\Box$, which can be solved with
$$t \mapsto \frac{4 (-331 + u^2)}{7 (-793 + u^2)}.$$
The corresponding specialization also has rank $5$.

But now we observe that in both conditions
 \begin{align*}
-(-4+7 t) (4+17 t)&=\Box \\
(-4 + 7 t) (-1324 + 5551 t)&=\Box,
\end{align*}
the same factor $(-4+7t)$ appears. This implies that both conditions can be satisfied simultaneously
since when we apply the solution of the first condition to the second,
we have to solve
$ -(-1863+539 u^2)=z^2$ in order to satisfy both conditions. But this genus $0$ curve
has rational points, i.e. $(u,z)=(13/7,2)$, so we can find the parametric solution
$$
u\mapsto  \frac{-7007-28 r+13 r^2}{7 (539+r^2)}.
$$


So we can satisfy both conditions with
\[
t\mapsto \frac{4 \left(3 r^2-14 r-5390\right) \left(10 r^2-14 r-1617\right)}{7 \left(72 r^4-182 r^3-13279 r^2+98098 r+20917512\right)}.
\]
By inserting this into $E$, we get the curve over $\Q(r)$ given by
\[ y^2=x^3+a_6(r) x^2 + b_6(r) x, \]
where
{\footnotesize
 \begin{align*}
a_6(r)&=3 (6637977907200 r^{16}-327957190299648 r^{15}-132939477324670464 r^{14}+\\&
1334557851651990784 r^{13}+73205200037549219248
   r^{12}-\\&
   1718125119359074284768 r^{11}-193538301177692188691736 r^{10}+\\&1905189555626165277886872 r^9+96624855648992854220247819
   r^8-\\&1026897170482503084781024008 r^7-56226940796444312350911834456 r^6+\\&269042619584910197333660344992 r^5+6178698341397939354226782536368
   r^4-\\&60712935351132451806093801142016 r^3-3259766993714464579957766495983104 r^2+\\&4334495070152077221968455683796992
   r+47287453161693896431461711200563200),\\
   b_6(r)&=7225344 (r-224)^2 (r+154)^2 (2 r-7)^2 (32 r+77)^2 \left(18 r^2+14 r+16709\right)^2 \\  &\left(24 r^2-1001 r+14014\right)^2 \left(26 r^2+1001
   r+12936\right)^2 \left(31 r^2-14 r+9702\right)^2\\& \left(72 r^4-182 r^3-13279 r^2+98098 r+20917512\right)^2.
\end{align*}
}%

The $x$-coordinates of six independent points of infinite order are
{\footnotesize
 \begin{align*}
x_1&=-53067 (r-224) (r+154) (2 r-7) (32 r+77) (24 r^2-1001 r+14014) \times\\&
(26 r^2+1001 r+12936) (72 r^4-182 r^3-13279 r^2+98098
   r+20917512)^2,\\
   x_2&=-48 (r-224) (r+154) (2 r-7) (32 r+77) (18 r^2+14 r+16709) \times\\&(31 r^2-14 r+9702) (2424 r^4-12922 r^3-3840473 r^2+6964958
   r+704222904)^2,\\
   x_3&=144 (16709 + 14 r + 18 r^2) (14014 - 1001 r + 24 r^2) (12936 +
   1001 r + 26 r^2) \times\\&(9702 - 14 r + 31 r^2) (155719256 - 490490 r +
   1032283 r^2 + 910 r^3 + 536 r^4)^2,\\
   x_4&=576 (16709 + 14 r + 18 r^2) (14014 - 1001 r + 24 r^2) (12936 +
   1001 r + 26 r^2)\times\\& (9702 - 14 r + 31 r^2) (434619416 + 2648646 r -
   841477 r^2 - 4914 r^3 + 1496 r^4)^2,\\
   x_5&=(12936 + 1001 r + 26 r^2)^2 (20917512 + 98098 r - 13279 r^2 -
   182 r^3 + 72 r^4)^2\times\\&
    \frac {88510464 (539 + r^2)^2 (-7007 - 28 r + 13 r^2)^2 (14014 -
     1001 r + 24 r^2)^2} {(285872664 + 2256254 r - 1029833 r^2 -
    4186 r^3 + 984 r^4)^2},\\
    x_6&=(9702 - 14 r + 31 r^2)^2 (20917512 + 98098 r - 13279 r^2 - 182 r^3 +
   72 r^4)^2\times\\&\frac {260112384 (539 + r^2)^2 (-539 + 1001 r + r^2)^2 (16709 +
     14 r + 18 r^2)^2} {(676332888 + 2256254 r + 418999 r^2 -
    4186 r^3 + 2328 r^4)^2}.
    \end{align*}
}%

Since the specification map is a homomorphism, to show that these six points are
independent, it suffices to find a rational number $r$ for which the points are specialized
to six independent $\mathbb{Q}$-rational points. It is easy to check that we may take, e.g. $r=1$.

Moreover, for $r=13$, the conditions of the Gusi\'c-Tadi\'c algorithm for finding injective specializations
(see \cite[Theorem 1.3]{GT2} and \cite{GT1}) are satisfied,
so the rank over $\mathbb{Q}(r)$ is exactly $6$.

This specialization also indicates that the points $P_1,P_2,P_3,P_4,P_5,P_6$, where $x(P_i)=X_i$ for $i=1,2,\ldots,6$,
do not generate the full Mordell-Weil group over $\mathbb{Q}(r)$ modulo torsion, but a subgroup of order $4$
(since the quotient of the corresponding regulators is $16$). Indeed, it holds that
$P_3+P_4+P_5 = 2R_5$ and $P_3+P_4+P_6=2R_6$, where

{\footnotesize
 \begin{align*}
x(R_5)&=1344(72r^4-182r^3-13279r^2+98098r+20917512)(2r-7)^2 (32r+77)^2 \times\\
& (18r^2+14r+16709)^2 (11r^2+14r+12936)^2 ,\\
x(R_6)&=1344(18r^2+14r+16709)(26r^2+1001r+12936)(31r^2-14r+9702) \times\\
&(24r^2-1001r+14014) (72r^4-182r^3-13279r^2+98098r+20917512) \times\\
&\frac{(r+154)^2 (32r+77)^2 (6r^2-91r+3332)^2}{(30r^2-1001r+17248)^2},
  \end{align*}
}%
and now the points $P_1,P_2,P_3,P_4,R_5,R_6$ generate the full Mordell-Weil group over $\mathbb{Q}(r)$ modulo torsion.


\section{Concluding remarks}

Let us mention that the heuristic of Park, Poonen, Voight and Wood
(see \cite[Section 8.3]{PPVW}) predicts that there are only finitely many elliptic curves over $\Q$
with torsion group $\Z/4\Z$ and rank greater than $7$, which indicates that our result
might be close to the best possible.

However, the recent paper by Dujella, Kazalicki and Peral \cite{DKP}
might indicate that the heuristic in \cite{PPVW} needs some adjustments, at least in
the case of curves with certain torsion groups.
Namely, they examined families with rank $3$ for torsion groups $\Z/8\Z$ and $\Z/2\Z \times \Z/6\Z$
and the experiments suggest that in these families, root
numbers $1$ and $-1$ (which conjecturally means even and odd ranks) are equidistributed,
which would imply that there are infinitely many curves with rank $\geq 4$ for these torsion
groups, while the heuristic in \cite{PPVW} predicts that there are only finitely many elliptic curves over
$\mathbb{Q}$ with torsion group $\Z/8\Z$ or $\Z/2\Z \times \Z/6\Z$ and rank greater than $3$.

\medskip

{\bf Acknowledgements.}
A.~D.~ acknowledges support from the QuantiXLie Center of Excellence, a project
co-financed by the Croatian Government and European Union through the
European Regional Development Fund - the Competitiveness and Cohesion
Operational Programme (Grant \linebreak KK.01.1.1.01.0004).




\bigskip
\bigskip

\begin{center}
{\bf Elipti\v{c}ka krivulja nad $\mathbb{Q}(u)$ s torzijom $\mathbb{Z}/ 4\mathbb{Z}$ i rangom $6$}
\end{center}

\bigskip

\begin{center}
{\it Andrej Dujella i Juan Carlos Peral}
\end{center}

\bigskip

\begin{center}
\begin{minipage}[c]{9.2cm}
{\small \hspace*{0.5cm} {\sc Sa\v{z}etak.}
U ovom radu prikazujemo konstrukciju
elipti\v{c}ke krivulje nad $\mathbb{Q}(u)$ s torzijskom grupom $\mathbb{Z}/ 4\mathbb{Z}$ i rangom $6$.
Prethodno su bili poznati samo primjeri ranga $5$ za takve krivulje.}%
\end{minipage}
\end{center}

\endarticle

\end{document}